\documentclass[12pt]{article}
\usepackage{mathrsfs}
\usepackage{amssymb}
\usepackage{amsmath}
\usepackage{amsfonts}
\usepackage{ulem}
\usepackage{graphicx,amsmath,amsfonts,amssymb,color,mathrsfs, amsmath,amsthm}
\usepackage{epsfig}
\newcommand{\E}{{\mathbb E}}

\newcommand{\chuhao}{\fontsize{19pt}{\baselineskip}\selectfont}

\newtheorem{thm}{Theorem}[section]
\newtheorem{prop}[thm]{Proposition}
\newtheorem{cor}[thm]{Corollary}
\newtheorem{lem}[thm]{Lemma}
\newtheorem{defi}[thm]{Definition}
\newtheorem{remark}[thm]{Remark}
\newtheorem{example}[thm]{Example}
\newtheorem{pb}[thm]{Problem}

\newcommand{\be}{\begin{eqnarray*}}
\newcommand{\ee}{\end{eqnarray*}}
\newcommand{\beq}{\begin{equation}}
\newcommand{\eeq}{\end{equation}}

\numberwithin{equation}{section}

\title{\bf\color{black} \chuhao{Weak Orlicz-Hardy Martingale Spaces}}
\author{Yong JIAO and Lian WU\\ \small{ (Institute of Probability and
Statistics, Central South University,
            Changsha
            410075, China)} \\ }
\date{\small Revised version, April 12, 2013}

\begin{document}
 \maketitle

 %%%%%%%%%%%%%%%%%%%%%%%%%%%%%%%%%%%%%%%%%%%%%%%%%%%%%%%%%%%%%%%%
%%%%%%%%%%%%%%%%%%%%%%%%%%%%%%%%%%%%%%%%%%%%%%%%%%%%%%%%%%%%%%%%

 \makeatletter
 \renewcommand{\@makefntext}[1]{#1}
 \makeatother \footnotetext{\noindent
 Partially supported by the National Natural Science Foundation of
China(11001273, 11150110456), the Research Fund for the Doctoral
Program of Higher Education of China (20100162120035) and
Postdoctoral Science Foundation of China and Central South
University.
}

%%%%%%%%%%%%%%%%%%%%%%%%%%%%%%%%%%%%%%%%%%%%%%%%%%%%%%%%%%%%%%%%
%%%%%%%%%%%%%%%%%%%%%%%%%%%%%%%%%%%%%%%%%%%%%%%%%%%%%%%%%%%%%%%%
%%%%%%%%%%%%%%%%%%%%%%%%%%%%%%%%%%%%%%%%%%%%%%%%%%%%%%%%%%%%%%%%
%%%%%%%%%%%%%%%%%%%%%%%%%%%%%%%%%%%%%%%%%%%%%%%%%%%%%%%%%%%%%%%%

 \begin{abstract}
In this paper, several weak Orlicz-Hardy martingale spaces
associated with concave functions are introduced, and some weak
atomic decomposition theorems for them are established. With the
help of weak atomic decompositions, a sufficient condition for a
sublinear operator defined on the weak  Orlicz-Hardy martingale
spaces to be bounded is given. Further, we investigate the duality
of weak Orlicz-Hardy martingale spaces and obtain a new
John-Nirenberg type inequality when the stochastic basis is regular.
These results can be regarded as weak versions of the Orlicz-Hardy
martingale spaces due to Miyamoto, Nakai and Sadasue.

\end{abstract}

%%%%%%%%%%%%%%%%%%%%%%%%%%%%%%%%%%%%%%%%%%%%%%%%%%%%%%%%%%%%%%%%
%%%%%%%%%%%%%%%%%%%%%%%%%%%%%%%%%%%%%%%%%%%%%%%%%%%%%%%%%%%%%%%%
\section{Introduction}
The Lebesgue's theory of integration has taken a center role in
modern analysis, which leads the more extensive classes of function
spaces and martingale spaces to naturally arise.  As is well known,
as a generalization of $L_p$-space, the Orlicz space was introduced
in \cite{orlicz}. Since then, the Orlicz spaces have been widely
used in probability, partial differential equations and harmonic
analysis; see \cite{arai,byun,or,rao,s}, and so forth. In
particular, Takashi, Eiichi and Gaku very recently studied the
Orlicz-Hardy martingale spaces in \cite{mns}, and using the atomic
decomposition they obtained some very interesting martingale
inequalities as well as the dulity, and proved a generalized
John-Nirenberg type inequality for martingale when the stochastic
basis is regular. Let us briefly recall the main results of
\cite{mns}.

Let $\mathcal {G}$ be the set of all functions
$\Phi:[0,\infty)\rightarrow [0,\infty)$ satisfying $\Phi(0)=0$,
 $\lim_{r\rightarrow\infty}\Phi(r)=\infty$. The
 Orlicz space $L_\Phi$ is defined as the collection of all measurable functions $f$
 with respect to $(\Omega,\mathcal{F},P)$ such that
 $\E(\Phi(c|f|))<\infty$ for some $c>0$ and
 $$\|f\|_{L_\Phi}=\inf\Big\{c>0:\E(\Phi(|f|/c))\leq 1 \Big\},$$
where $\E$ denotes the expectation with respect to $\mathcal{F}$.
For $q\in [1,\infty)$ and a function $\phi:(0, \infty)\rightarrow(0,
\infty)$, the generalized Campanato martingale spaces $\mathcal
{L}_{q,\phi}$ is defined by $\mathcal {L}_{q,\phi}= \{f\in L_{q}:
\|f\|_{\mathcal {L}_{q,\phi}} < \infty \}$, where
\begin{equation}\label{1}\|f\|_{\mathcal {L}_{q,\phi}} = \sup_{n\geq 1}\sup_{A\in\mathcal
{F}_n}\frac{1}{\phi(P(A))}\Bigg(\frac{1}{P(A)}\int_A|f-\E_{n} f|^q
dP\Bigg)^{1/q},\end{equation} with the convention that $\E_0f=0.$ We
refer to the recent paper \cite{ns} for the Morrey-Campanato spaces.
Denote by $\mathcal{G}_{\ell}$ the set
$$\mathcal{G}_{\ell}=\Big\{\Phi\in \mathcal{G}: \exists\, c_\Phi\geq 1\,
{\rm and}\, \ell\in(0,1]\,\,{\rm s.t. }\Phi(tr)\leq
c_\Phi\max\{t^{\ell},t\}\Phi(r) \,\,{\rm for}\,t,r \in
[0,\infty)\Big\}.$$ Then for $\Phi\in \mathcal{G}_{\ell}$ and
$\phi(r)=\frac{1}{r\Phi^{-1}(1/r)}$, where and in what follows
$\Phi^{-1}$ denotes the inverse function of $\Phi,$ the following
duality holds,
$$\big(H_\Phi^{s}\big)^*=\mathcal {L}_{2,\phi}\,.$$ See Section 2 for
the notation $H_\Phi^{s}.$ Moreover, the John-Nirenberg type
inequality holds when the stochastic basis is regular, namely,
$\mathcal {L}_{q,\phi}$ are equivalent for all $1\leq q<\infty.$ It
should be mentioned that Miyamoto, Nakai and Sadasue's results above
are exactly the generalization in \cite{weisz2} due to Weisz when
$\Phi(t)=t^p$, $0<p\leq1$.

The main goal of this present paper is to deal with the weak
Orlicz-Hardy martingale spaces, which are more inclusive classes
than Orlicz-Hardy martingale spaces, and give the weak version of
Miyamoto, Nakai and Sadasue's results. In 2008 the weak Orlicz
spaces and weak Orlicz-Hardy martingale spaces generated by nice
young functions satisfying the $M_\triangle$-condition were first
introduced in \cite{jiao1}, and interpolation theorems and
inequalities were proved for these spaces; Liu et al investigated
the boundedness of some sublinear operators defined on weak
Orlicz-Hardy martingale spaces in \cite{liu1} and the first named
author studied some embedding relationships between them in
\cite{jiao3} in 2011; however, the existing results about weak
Orlicz-Hardy martingale spaces are all associated with convex
functions. In the present paper we are interested in the case $\Phi$
is not convex. We denote
$$t_{\phi}^{q}(x)=\frac{1}{\phi(x)}x^{-1/q}\sup_{P(\nu<\infty)\leq x}\|f-f^{\nu}\|_q,$$
where $\nu$ is a stopping time and $f^{\nu}$ is the stopped
martingale. Very differently from \eqref{1}, we define the weak
generalized Campanato martingale space $w\mathcal {L}_{q,\phi}$ as
follows.
\begin{defi} For $q\in[1,\infty)$ and a function
$\phi:(0, \infty)\rightarrow(0, \infty)$, let
$$w\mathcal {L}_{q,\phi}=\Big\{f\in L_{q}:\|f\|_{w\mathcal {L}_{q,\phi}}
=\int_0^{\infty}\frac{t_{\phi}^q (x)}{x} dx <\infty \Big\}.$$
\end{defi}
Then for $\Phi\in \mathcal{G}_{\ell}$ and
$\phi(r)=\frac{1}{r\Phi^{-1}(1/r)}$, we have
$$\big(w\mathscr{H}_\Phi^{s}\big)^*=w\mathcal {L}_{2,\phi}\,.$$ See Section 2 for
the notation of $w\mathscr{H}_\Phi^{s}.$ Furthermore, we by the
duality obtain the weak type John-Nirenberg inequality when the
stochastic basis is regular. That is, $w\mathcal {L}_{q,\phi}$ are
equivalent for all $1\leq q<\infty.$ We note that our theorems can
deduce Weisz's main results in \cite{weisz1}.

It is well known that the method of atomic decompositions plays an
important role in martingale theory; see for example,
 \cite{hou,jiao2,weisz3,weisz4}. In the present paper, the
important step is to establish the weak atomic characterizations of
weak Orlicz-Hardy martingale spaces. To this end, the main
difficulty encountered is that there is no similar result to replace
Lemma 3.1 or Remark 3.2 in \cite{mns}. Inspired by \cite{weisz1}, we
adopt a different method and apply some technical estimates.
Particularly, we note that $\frac{\Phi^{-1}(t)}{t^{p}}$ are
increasing and $\frac{\Phi^{-1}(t)}{t^{q}}$ are decreasing on
$(0,\infty)$ for $\Phi\in\mathcal{G}_{\ell}$ with $q<\infty$, where
$p=p_{\Phi^{-1}}$ and $q=q_{\Phi^{-1}}$ denote the lower index and
the upper index of convex function $\Phi^{-1}$,  respectively; see
also Section 2 for the definitions of the lower and upper index.

This paper is organized as follows. Section 2 is on preliminaries
and notations. Section 3 is devoted to the weak atomic
decompositions of weak Hardy-Orlicz martingale spaces. By the atomic
decompositions, a sufficient condition for a sublinear operator
defined on weak Hardy-Orlicz martingale spaces to be bounded is
given in Section 4. In Section 5, we deduce the new John-Nirenberg
type inequality by duality.

We end this section by an open question. For $q\in [1,\infty)$ and a
function $\phi:(0, \infty)\rightarrow(0, \infty)$, define $
{L}_{q,\phi}= \{f\in L_{q}: \|f\|_{ {L}_{q,\phi}} < \infty \}$,
where $$\|f\|_{{L}_{q,\phi}} =
\sup_{\nu}\frac{1}{\phi(P(\nu<\infty))}
\Bigg(\frac{1}{P(\nu<\infty)}\int_{\{\nu<\infty\}}|f- f^{\nu}|^q
dP\Bigg)^{1/q}.$$  The supremum is taken all the stopping times
$\nu.$ Then by the Proposition 2.12 in \cite{mns},
\begin{equation}\label{2} \|f\|_{\mathcal{L}_{q,\phi}}\leq\|f\|_{{L}_{q,\phi}}\leq C_\Phi\|f\|_{\mathcal{L}_{q,\phi}},\end{equation}
where $\phi$ satisfies $\phi(r)\leq c_\Phi\phi(s)$ for $0<r\leq
s<\infty.$ Now we denote
$$u_{\phi}^{q}(x)=\frac{1}{\phi(x)}x^{-1/q}\sup_{n\geq 1}\sup_{A\in\mathcal
{F}_n}\sup_{P(A)\leq x}\Bigg(\int_A|f-\E_{n} f|^q dP\Bigg)^{1/q},$$
where $\nu$ is a stopping time and $f^{\nu}$ is the stopped
martingale. It is natural to define
$$w {L}_{q,\phi}=\Big\{f\in L_{q}:\|f\|_{w{L}_{q,\phi}}
=\int_0^{\infty}\frac{u_{\phi}^q (x)}{x} dx <\infty \Big\}.$$ In the
time of writing this paper, we do not know if there is a result
similar to \eqref{2}.

Throughout this paper, $\textbf{Z}$ and $\textbf{N}$ denote the
integer set and nonnegative integer set, respectively. We denote by
$C$ the positive constant, which can vary from line to line.

\section{Preliminaries}
 Let $(\Omega,\mathcal{F},P)$ be a probability space, and $\{\mathcal {F}_n\}_{n\geq0}$ be a non-decreasing sequence of
 sub-$\sigma $-algebras of $\mathcal {F}$ such that $\mathcal
 {F}=\mathcal {\sigma }(\bigcup_{n\geq 0}\mathcal {F}_n)$. The
 expectation operator and the conditioned expectation operator are
 denoted by $\E$ and $\E_n$, resp. For a martingale $f=(f_n)_{n\geq 0}$
 relative to $(\Omega,\mathcal {F},P;(\mathcal {F}_n)_{n\geq 0})$, we
 denote its martingale difference by $df_i=f_i-f_{i-1}$ ( $i\geq 0$, with convention
 $f_{-1}=0$). Then the maximal function, the quadratic variation and
 the conditional quadratic variation of martingale $f$ are defined
 by $$M_n(f) =\sup_{0\leq i \leq n}|f_i|,\qquad M(f)=\sup_{i\geq 0}|f_i|,$$
 $$S_n(f)=\Bigg(\sum_{i=0}^n |df_i|^2\Bigg)^{1/2},\qquad S(f)=\Bigg(\sum_{i=0}^\infty |df_i|^2\Bigg)^{1/2},$$
 $$s_n(f)=\Bigg(\sum_{i=0}^n \E_{i-1}|df_i|^2\Bigg)^{1/2},\qquad s(f)=\Bigg(\sum_{i=0}^\infty \E_{i-1}|df_i|^2\Bigg)^{1/2}.$$
 The stochastic basis $(\mathcal {F}_n)_{n\geq 0}$ is said to be
 regular, if for $n\geq 0$ and $A\in \mathcal {F}_n$, there exists
 $B\in \mathcal {F}_{n-1} $ such that $A\subset B$ and $P(B)\leq R
 P(A)$, where $R$ is a positive constant independent of $n$. A martingale is
 said to be regular if it is adapted to a regular $\sigma$-algebra
 sequence. This amounts to saying that there exists a constant $R>0$
 such that $$f_n\leq Rf_{n-1}$$  for all nonnegative
 martingales $(f_n)_{n\geq0}$ adapted to the stochastic basis $(\mathcal {F}_n)_{n\geq
 0}$.

  Recall that $\mathcal {G}$ is the collection of all functions $\Phi:[0,\infty)\rightarrow [0,\infty)$
  satisfying $\Phi(0)=0$,
 $\lim_{r\rightarrow\infty}\Phi(r)=\infty$. A function $\Phi$ is said to satisfy $\Delta_2$-condition,
 denoted by $\Phi\in \Delta_2$, if there exists a constant $C>0$ such that $\Phi(2t)\leq
C\Phi(t), \forall t>
 0.$ A function $\Phi:[0,\infty)\rightarrow [0,\infty)$
 is said to be subadditive if $\Phi(r+s)\leq \Phi(r)+\Phi(s)$ for all $r,s\in
 [0,\infty)$. We note that all concave functions are subadditive.
Let $\Phi_1$ and $\Phi_2$ belong to $\mathcal {G}$, which are said
to be equivalent if there exists a constant $C\geq1$ such that
$\Phi_1(t)/C\leq \Phi_2(t)\leq C\Phi_1(t)$ for all $t\geq0.$
\begin{defi}Let $\Phi\in\mathcal {G}$, then the weak
 Orlicz space $wL_\Phi$ is defined as the set of all measurable functions $f$
 with respect to $(\Omega,\mathcal{F},P)$ such that
 $\|f\|_{wL_\Phi}<\infty,$ where
 $$\|f\|_{wL_\Phi}:=\inf\Big\{c>0:\sup_{t>0}\Phi\Big(\frac{t}{c}\Big)P(|f|>t)\leq 1 \Big\}.$$
\end{defi}
If $\Phi(t)=t^p$, $0<p<\infty$, then $wL_\Phi=wL_p$, where the weak
$L_p$ space $wL_p$ consists of all measurable functions $f$ for
which
$$\|f\|_{wL_p}:=\sup_{t>0}t P(|f|>t)^{1/p}<
\infty.$$ It was proved in \cite{liu1} that the functional
$\|\cdot\|_{wL_\Phi}$ is a
 complete quasi-norm on $wL_\Phi$ when $\Phi\in\Delta_2$ and $\Phi$ is convex.
In this paper we are interested in
 the case $\Phi$ is not convex. We assume that $\Phi$ is of lower
 type $\ell$ for some $\ell\in(0,1]$ and upper type 1, i.e., there exist
 a constant $c_\Phi\in[1,\infty)$ and some $\ell\in(0,1]$ such that
 $$\Phi( tr)\leq c_\Phi\max\{t^\ell,t\}\Phi(r)\quad {\rm for }\quad t,r\in[0,\infty).$$
 Let $\mathcal{G}_\ell$ be the set of all $\Phi\in\mathcal{G}$
 satisfying the above inequality. For example,
 $\Phi(t)=t^p\big(\log(e+t)\big)^q$ is in $\mathcal{G}_\ell$ if
 $0<\ell\leq p<1$ and $q\geq0.$
 Let $\Phi\in \mathcal {G}_\ell $, from
 \cite{peetre} we know that $\Phi$ is equivalent to a concave
 function in $\mathcal {G}_\ell $. Further, we can verify the functionals $\|\cdot\|_{wL_\Phi}$ and
 $\|\cdot\|_{wL_\Psi}$ are equivalent if $\Phi,\Psi\in \mathcal
 {G}_\ell$ are equivalent. Therefore, we always assume that $\Phi\in\mathcal
 {G}_\ell$ is concave in our theorems below. Thus $\Phi$ is subadditive, increasing, continuous and bijective from
 $[0,\infty)$ to itself when $\Phi\in\mathcal
 {G}_\ell$.

Obviously, $\mathcal {G}_\ell \subset \Delta_2$. It was shown in
\cite{mns} that for
 any concave function $\Phi\in\mathcal
 {G}_\ell$, $L_\Phi$ is a $\ell$-quasi norm. Here the functional $\|\cdot\|_{wL_\Phi}$ satisfies the
 following properties:

 {\rm(i)} $\|f\|_{wL_\Phi}\geq 0$, and $\|f\|_{wL_\Phi}=0$ if and only if $f=0$;

 {\rm(ii)} $\|cf\|_{wL_\Phi}=|c|\|f\|_{wL_\Phi}$;

 {\rm(iii)} $\|f+g\|_{wL_\Phi}\leq
 (2c_\Phi)^{1/\ell}\cdot2(\|f\|_{wL_\Phi}+\|g\|_{wL_\Phi})$.

\noindent Indeed, (i) and (ii) can be proved easily. Here we only
prove the
 generalized triangle inequality, namely (iii). Suppose that
 $\|f\|_{wL_\Phi}=a$, $\|g\|_{wL_\Phi}=b$, $a,b>0$, and $K=(2c_\Phi)^{1/\ell}$. Then $\forall
 t>0$,
 \be
 \Phi\Big(\frac{t}{K\cdot2(a+b)}\Big)P(|f+g|>t) & \leq
 & \Phi\Big(\frac{t}{K\cdot2(a+b)}\Big)\Bigg(P\Big(|f|>\frac{t}{2}\Big)+P\Big(|g|>\frac{t}{2}\Big)\Bigg)\\
 & \leq & c_\Phi \cdot
 \frac{1}{2c_\Phi}\Bigg(\Phi\Big(\frac{t}{2a}\Big)P\Big(|f|>\frac{t}{2}\Big)+\Phi\Big(\frac{t}{2b}\Big)P\Big(|g|>\frac{t}{2}\Big)\Bigg)\\
 & \leq & 1.
 \ee

\begin{prop}
 Let $\Phi \in \mathcal {G}_{\ell}$, then $L_1\subset L_\Phi\subset wL_\Phi$.
 \end{prop}

 \noindent {\bf Proof.} Suppose that $f\in L_1$, then
 \be
 \int_{\Omega}\Phi\Big(\frac{|f|}{\|f\|_1}\Big)d P
 & = &\int_{\{|f|\leq\|f\|_1\}}\Phi\Big(\frac{|f|}{\|f\|_1}\Big)d p + \int_{\{|f|>\|f\|_1\}}\Phi\Big(\frac{|f|}{\|f\|_1}\Big)d p\\
 & \leq & \int_{\{|f|\leq\|f\|_1\}}\Phi(1)d p + \int_{\{|f|>\|f\|_1\}}c_\Phi \max\Big\{\frac{|f|}{\|f\|_1}, \Big(\frac{|f|}{\|f\|_1}\Big)^\ell \Big \}\Phi(1)d p\\
 & \leq & \Phi(1)+c_\Phi \Phi(1) \int_{\{|f|>\|f\|_1\}} \frac{|f|}{\|f\|_1}d
 P\\
 & \leq & \Phi(1) + c_\Phi \Phi(1).
 \ee
 Taking $C_0=\max\{\Phi(1)+c_\Phi \Phi(1), 1\}$, then
 \be
 \int_{\Omega}\Phi\Big(\frac{|f|}{(C_0\cdot c_\Phi)^{1/\ell}\|f\|_1}\Big)d
 P & \leq & c_\Phi\cdot\frac{1}{C_0c_\Phi}\int_{\Omega}\Phi\Big(\frac{|f|}{\|f\|_1}\Big)d
 P\leq  1,
 \ee
 which means
 $$\|f\|_{L_\Phi}\leq (C_0\cdot c_\Phi)^{1/\ell}\|f\|_1.$$
Suppose that $f\in L_\Phi$ and $t>0$. Then \be
 \Phi\Big(\frac{t}{\|f\|_{L_\Phi}}\Big)P(|f|>t)
 & \leq & \int_{\{|f|>t\}}\Phi\Big(\frac{|f|}{\|f\|_{L_\Phi}}\Big)d P
  \leq  \int_{\Omega}\Phi\Big(\frac{|f|}{\|f\|_{L_\Phi}}\Big)d P\leq
  1,\ee  which implies $L_\Phi\subset wL_\Phi$. The proof is complete.

\begin{remark} It was proved in \cite{rao} that the Orlicz space $L_\Phi$ has
absolute continuous norm when $\Phi\in\Delta_2,$ that is,
$$\lim_{P(A)\rightarrow0}\|f\chi_A\|_{L_\Phi}=0,\quad \quad \forall f\in L_\Phi.$$
But not every element in $wL_\Phi$ has absolute continuous quasi
norm in spire of $\Phi\in\Delta_2$. For instance, let $\Omega=(0,1]$
and $P$ be the Lebesgue measure on $\Omega$. Consider $wL_p$
($0<p<\infty$) and function $f(x)=x^{-1/q}$. A simple computation
shows that $f\in wL_p$ when $q\geq p$, and $f$ has absolutely
continuous norm in $wL_p$ when $q>p$, but it has not when $q=p$.
\end{remark}

\begin{defi}Let $w\mathscr{L}_\Phi$ be the set of all $f\in  wL_\Phi$ having the absolute continuous quasi norm defined by
 $$w\mathscr{L}_\Phi:=\Big\{f\in  wL_\Phi: \lim_{P(A)\rightarrow0}\|f\chi_A\|_{wL_\Phi}=0\Big\}.$$
\end{defi}

It is easy to check $w\mathscr{L}_\Phi$ is a linear space. Moreover
$w\mathscr{L}_\Phi$ is a closed subspace of $wL_\Phi$ when $\Phi \in
\mathcal {G}_{\ell}.$ Indeed, suppose that $f_n\in
w\mathscr{L}_\Phi$, and $\|f_n-f\|_{wL_\Phi}\rightarrow 0$ as
$n\rightarrow \infty$. Then $f\in wL_\Phi$ and
 \be
 \|f\chi_A\|_{wL_\Phi}& \leq &
 (2c_\Phi)^{1/\ell}\cdot2(\|f_n \chi_A\|_{wL_\Phi} + \|(f_n
 -f)\chi_A\|_{wL_\Phi} )\\
 & \leq & (2c_\Phi)^{1/\ell}\cdot2(\|f_n \chi_A\|_{wL_\Phi} + \|f_n
 -f\|_{wL_\Phi}).
 \ee
Since $f_n\in w\mathscr{L}_\Phi$ has absolute continuous norm and
$\|f_n-f\|_{wL_\Phi}\rightarrow 0$ as $n\rightarrow \infty$, we get
$$ \lim_{P(A)\rightarrow 0}\|f\chi_A\|_{wL_\Phi}= 0,$$ namely, $f\in
w\mathscr{L}_\Phi$, which means $w\mathscr{L}_\Phi$ is closed.
Moreover, $L_1 \subset L_\Phi \subset w\mathscr{L}_\Phi$ when $\Phi
\in \mathcal {G}_{\ell}.$

\smallskip
 The following is an extension of Lebesgue  controlled convergence
 theorem, which will be used to describe the quasi norm convergence
 (See Remark 3.2).
\begin{prop}(See \cite{liu2}, Theorem 3.2)
 Let $\Phi \in \mathcal {G}_{\ell}$. $f_n, g\in w\mathscr{L}_\Phi$ and $|f_n|\leq
 g$. If $f_n$ converges to $f$ almost everywhere, then $$\lim_{n\rightarrow\infty}\|f_n-f\|_{w\mathscr{L}_\Phi}=0.$$
 \end{prop}

 \noindent In order to describe our results, we need the lower index and upper
 index of $\Phi$. Let $\Phi\in\mathcal {G}$, the lower index and the upper
 index of $\Phi$ are respectively defined by
 $$p_\Phi=\inf_{t>0}\frac{t\Phi^{'}(t)} {\Phi(t)},\qquad q_\Phi=\sup_{t>0}\frac{t\Phi^{'}(t)} {\Phi(t)}.$$
It is well known that $1\leq p_\Phi\leq q_\Phi\leq\infty$ if $\Phi$
is convex and  $0< p_\Phi\leq q_\Phi\leq1$ if $\Phi$ is concave.

\smallskip
We note that the important observation below.
\begin{lem} Let $\Phi$ be concave and $q_{\Phi^{-1}}<\infty$. Denote $p=p_{\Phi^{-1}}$,
$q=q_{\Phi^{-1}}$. Then $\frac{\Phi^{-1}(t)}{t^p}$,
$\frac{\Phi(t)}{t^{1/q}}$ are increasing on $(0,\infty)$ and
$\frac{\Phi^{-1}(t)}{t^q}$, $\frac{\Phi(t)}{t^{1/p}}$ are decreasing
on $(0,\infty)$.
\end{lem}

\noindent {\bf Proof.} It is easy to see that $\Phi^{-1}$ is convex.
Thus $1\leq p \leq q <\infty$. From \cite{rao}, we obtain that
$\frac{\Phi^{-1}(t)}{t^p}$ is increasing on $(0,\infty)$ and
$\frac{\Phi^{-1}(t)}{t^q}$ is decreasing on $(0,\infty)$. Replacing
$t$ with $\Phi(t)$, we immediately get that
$\frac{\Phi(t)}{t^{1/q}}$ is increasing on $(0,\infty)$ and
$\frac{\Phi(t)}{t^{1/p}}$ is decreasing on $(0,\infty)$.

\bigskip

 We now introduce the weak Orlicz-Hardy martingale spaces. Denote by $\Lambda$ the collection of all sequences
 $(\lambda_n)_{n\geq0}$ of non-decreasing, non-negative and adapted
 functions with $\lambda_\infty=\lim_{n\rightarrow
 \infty}\lambda_n$. As usual, the weak Orlicz-Hardy martingale
 spaces are defined as follows:
 $$wH_\Phi = \{f=(f_n)_{n\geq 0}:\|f\|_{wH_\Phi}=\|M(f)\|_{wL_\Phi}<\infty\};$$
 $$wH_\Phi^{S} = \{f=(f_n)_{n\geq
 0}:\|f\|_{wH_\Phi^{S}}=\|S(f)\|_{wL_\Phi}<\infty\};$$
 $$ wH_\Phi^{s} = \{f=(f_n)_{n\geq 0}:\|f\|_{wH_\Phi^{s}}=\|s(f)\|_{wL_\Phi}<\infty\};$$
 $$w\mathcal{Q}_\Phi = \{f=(f_n)_{n\geq 0}:\exists(\lambda_n)_{n\geq0}\in\Lambda,\,\, \textrm{s.t.}\,\, S_n(f)
 \leq \lambda_{n-1}, \lambda_\infty\in wL_\Phi\},$$
 $$\|f\|_{w\mathcal{Q}_\Phi}=\inf_{(\lambda_n) \in \Lambda}\|\lambda_\infty\|_{wL_\Phi};$$
 $$w\mathcal{D}_\Phi = \{f=(f_n)_{n\geq 0}:\exists(\lambda_n)_{n\geq0}\in\Lambda,\,\, \textrm{s.t.}\,\,
 |f_n| \leq \lambda_{n-1}, \lambda_\infty\in wL_\Phi\},$$
 $$\|f\|_{w\mathcal{D}_\Phi}=\inf_{(\lambda_n) \in \Lambda}\|\lambda_\infty\|_{wL_\Phi}.$$

 \begin{remark} We can get the Orlicz-Hardy martingale space $H_\Phi^{s}$ when
 $\|s(f)\|_{wL_\Phi}$ is replaced by $\|s(f)\|_{L_\Phi}$ in the
 definition above. In order to describe the duality, we define
 $$w\mathscr{H}_\Phi^s =  \{f=(f_n)_{n\geq 0}:s(f)\in w\mathscr{L}_\Phi\}.$$
It is easy to see $w\mathscr{H}_\Phi^s$ is a closed subspace of
$wH_\Phi^s$.
 Similarly, we have  $w\mathscr{H}_\Phi$ and
 $w\mathscr{H}_\Phi^{S}$, which are closed
 subspaces of $wH_\Phi$ and $wH_\Phi^{S}$, respectively.
 \end{remark}

\begin{defi}A measurable function $a$ is said to be a
w-1-atom ( or w-2-atom, w-3-atom, resp.) if there exists a stopping
time $\nu$ such that

{\rm (a1)}\quad $a_n=E_n a =0$ if $\nu\geq n$,

{\rm (a2)}\quad $\|s(a)\|_\infty < \infty$ (or $\|S(a)\|_\infty <
\infty$, $\|M(a)\|_\infty < \infty$ resp.).
\end{defi}

\noindent  Now we define the weak Orlicz-Hardy spaces associated
with weak atoms .
\begin{defi}Let $\Phi \in \mathcal {G}_\ell$ with $\ell
\in (0,1]$. We define $wH_{\Phi,at}^s$ ($wH_{\Phi,at}^S$,
$wH_{\Phi,at}$ resp.) as the space of all $f\in wL_\Phi$ which admit
a decomposition
 $$f=\sum_{k\in \textbf{Z}} a^k\quad \quad a.e.$$ with for each $k\in\textbf{Z}$,
 $a^k$ is a w-1-atom (w-2-atom, w-3-atom resp.) and satisfying $s(a^k)$($S(a^k)$, $M(a^k)$ resp.) $\leq A\cdot 2^k)$ for some $A>0$, and
$$\inf \Big\{c>0:\sup_{k\in \textbf{Z}}\Phi\Big(\frac{2^k}{c}\Big)P(\nu_k < \infty)\leq
1\Big\}<\infty,$$

\noindent where $\nu_k$ is the stopping time corresponding to $a^k$.

\noindent Moreover, define
 $$\|f\|_{wH_{\Phi,at}^s}\quad (\|f\|_{wH_{\Phi,at}^S},\|f\|_{wH_{\Phi,at}} \textrm{resp}.)= \inf \inf \Big\{c>0:\sup_{k\in
\textbf{Z}}\Phi\Big(\frac{2^k}{c}\Big)P(\nu_k < \infty)\leq
1\Big\}$$ where the infimum is taken over all decompositions of $f$
described above.
\end{defi}
 Recall that for $q\in[1,\infty)$ and a function
$\phi:(0, \infty)\rightarrow(0, \infty)$,
$$w\mathcal {L}_{q,\phi}:=\Big\{f\in L_{q}:\|f\|_{w\mathcal {L}_{q,\phi}}
=\int_0^{\infty}\frac{t_{\phi}^q (x)}{x} dx <\infty \Big\},$$
 where
 $$t_{\phi}^{q}(x)=\frac{1}{\phi(x)}x^{-1/q}\sup_{P(\nu<\infty)\leq x}\|f-f^{\nu}\|_q,$$
 and $\nu$ is a stopping time. Then we have

 \begin{prop} If $1\leq q_1 \leq q_2 <\infty$, then
 $$\|f\|_{w\mathcal{L}_{q_1,\phi}}\leq \|f\|_{w\mathcal {L}_{q_2,\phi}}.$$
 \end{prop}

 \noindent {\bf Proof.} By H\"{o}lder's inequality,
 \be
 t_\phi^{q_1}(x) & = & \frac{1}{\phi(x)}x^{-1/q_1}\sup_{P(\nu<\infty)\leq
 x}\Big(\E(|f-f^\nu|^{q_1}\chi(\nu<\infty))\Big)^{1/q_1}\\
 & \leq & \frac{1}{\phi(x)}x^{-1/q_1}\sup_{P(\nu<\infty)\leq
 x}\big(\E(|f-f^\nu|^{q_2}\big)^{1/q_2}P(\nu<\infty)^{(1-q_1/q_2)(1/q_1)}\\
 & \leq & \frac{1}{\phi(x)}x^{-1/q_2} \sup_{P(\nu<\infty)\leq
 x}\|f-f^{\nu}\|_{q_2} = t_\phi^{q_2}(x),
 \ee
 which shows the proposition.

\section{Weak Atomic Decompositions} We now are in a position to prove the weak atomic decomposition of
the weak martingale Orlicz-Hardy spaces.

\begin{thm} Let $\Phi \in \mathcal {G}_\ell$ with $\ell
\in (0,1]$ and $q_{\Phi^{-1}}<\infty$.  Then $f \in wH_{\Phi}^{s}$
if and only if there exist a sequence of w-1-atoms $\{a^k\}_{k\in
\textbf{Z}}$ and corresponding stopping times $\{\nu_k\}_{k\in
\textbf{Z}}$ such that

{\rm(i)} $f_n = \sum_{k\in \textbf{Z}}E_n a^k\quad a.e.$, $\forall
n\in \textbf{N}$;

{\rm(ii)} $s(a^k)\leq A\cdot 2^{k}$ for some $A>0$ and
$$\inf \Big\{c>0:\sup_{k\in \textbf{Z}}\Phi\Big(\frac{2^k}{c}\Big)P(\nu_k < \infty)\leq
1\Big\}<\infty. $$

\noindent Moreover,
\begin{equation}\|f\|_{wH_{\Phi}^{s}}\approx \inf \inf \Big\{c>0:\sup_{k\in
\textbf{Z}}\Phi\Big(\frac{2^k}{c}\Big)P(\nu_k < \infty)\leq 1\Big\},
\end{equation} where the infimum is taken over
all the preceding decompositions of $f$. Consequently,
$$wH_\Phi^s = wH_{\Phi,at}^s \qquad with\;\,equivalent\;\,quasi\;\,norms.$$

\end{thm}

\noindent {\bf Proof.} Assume that $f =(f_n)_{n\geq 0} \in
wH_{\Phi}^{s}$. For $k\in\textbf{Z}$ and $n\geq 0 $, let
$$\nu_k = \inf\{n: s_{n+1}(f)>2^k\},\qquad a_n^k=f_n^{\nu_{k+1}}- f_n^{\nu_{k}}$$
Then it's clear that $\{\nu_k\}$ is nondecreasing and that for any
fixed $k\in \textbf{Z}$, $a^k=(a_n^k)_{n\geq 0}$ is a martingale.
Further we can see
$$s(f^{\nu_k})= s_{\nu_k}(f)\leq 2^k$$
and
$$\sum_{k\in \textbf{Z}}a_n^k=\sum_{k\in \textbf{Z}}(f_n^{\nu_{k+1}}-f_n^{\nu_{k}})=f_n,\;\; \textrm{a.e.}$$
for all $n\geq0$. Since $s(f^{\nu_k})\leq 2^k$, we have \be
s(a^k)&=&\Big(\sum_{n=1}^{\infty}\E_{n-1}|da_n^k|^2\Big)^{1/2}=\Big(\sum_{n=1}^{\infty}\E_{n-1}|df_n^{\nu_{k+1}}-df_n^{\nu_{k}}|^2\Big)^{1/2}\\
&=&\Big(\sum_{n=1}^{\infty}\E_{n-1}|df_n\chi_{\{\nu_k<n\leq\nu_{k+1}\}}|^2\Big)^{1/2}=\Big(\sum_{n=1}^{\infty}\chi_{\{\nu_k<n\leq\nu_{k+1}\}}\E_{n-1}|df_n|^2\Big)^{1/2}\\
&\leq& s(f^{\nu_{k+1}})\leq 2^{k+1}= 2\cdot 2^k. \ee Thus,
$(a_n^k)_{n\geq0}$ is $L_2$-bounded. Denote the limit still by
$a^k$. Then $a_n^k=\E_n a^k$ for all $n\geq 0$. For $\nu_k\geq n$,
$a_n^k = f_n^{\nu_{k+1}}- f_n^{\nu_k}=0.$ So $a^k$ is a w-1-atom and
(i) holds. Since $\{\nu_k<\infty\}=\{s(f)>2^k\}$, for any $k \in
\textbf{Z}$ we have
$$\Phi\Bigg(\frac{2^k}{\|f\|_{wH_{\Phi}^{s}}}\Bigg)P(\nu_k
<\infty)=\Phi\Bigg(\frac{2^k}{\|f\|_{wH_{\Phi}^{s}}}\Bigg)P(s(f)>2^k)\leq
1.$$ Thus $$\inf \Big\{c>0:\sup_{k\in
\textbf{Z}}\Phi\Big(\frac{2^k}{c}\Big)P(\nu_k < \infty)\leq
1\Big\}\leq \|f\|_{wH_{\Phi}^{s}}<\infty.$$

The main part of the proof is the converse. Suppose that there
exists a sequence of w-1-atoms such that (i) and (ii) hold. Let
$$M=\inf \Big\{c>0:\sup_{k\in \textbf{Z}}\Phi\Big(\frac{2^k}{c}\Big)P(\nu_k < \infty)\leq
1\Big\}.$$ Without loss of generality, we may assume that
$A=2^{N_A}$, where $N_A\geq0$ (Since for any $A>0$, there exist
$N_A\geq0$ such that $A\leq 2^{N_A}$). For any $\lambda>0$, choose
$j\in \textbf{Z}$ such that $2^j\leq\lambda<2^{j+1}$. Now let
\begin{equation}f_n=\sum_{k=-\infty}^{\infty}a_n^k=\sum_{k=-\infty}^{j-1}a_n^k+\sum_{k=j}^{\infty}a_n^k:=g_n+h_n\qquad(n\in
\textbf{N}).\end{equation} Then we have $s(f)\leq s(g)+s(h)$ by the
sublinearity of $s(f)$. And thus
$$P(s(f)>2A\lambda)\leq P(s(g)>A\lambda)+P(s(h)>A\lambda).$$
By (ii) we obtain
$$s(g)\leq \sum_{k=-\infty}^{j-1}s(a^k)\leq
\sum_{k=-\infty}^{j-1}A\cdot 2^k\leq A\cdot2^j\leq A\lambda.$$ So
$P(s(g)>A\lambda)=0$. Since $a_n^k=\E_n a^k = 0$ on the set
$\{\nu_k\geq n\}$, thus $s(a^k)= 0$ on $\{\nu_k =\infty\}$. Denote
$p=p_{\Phi^{-1}}$, $q=q_{\Phi^{-1}}$, then by Lemma 2.6 we have
 \be
 \Phi\Big(\frac{2A\lambda}{4
 M}\Big)P\big(s(f)>2A\lambda\big) & \leq &
 \Phi\Big(\frac{A\lambda}{M}\Big)P\big(s(h)>A\lambda\big)\\
 & \leq & \Phi\Big(\frac{A\lambda}{M}\Big)\sum_{k=j}^{\infty}P(\nu_k <\infty)\\
 & \leq & \Phi\Big(\frac{A\lambda}{M}\Big)\sum_{k=j}^{\infty}\frac{1}{\Phi
 \Big(\frac{2^k}{M}\Big)}.\ee
 Using the assumption
 $A=2^{N_A}$, we obtain $2^j\leq 2^{N_A + j}\leq A\lambda =
2^{N_A}\lambda <2^{N_A +j+1}$, where $N_A\geq 0$. Thus,
 \be
 \Phi\Big(\frac{2A\lambda}{4
 M}\Big)P(s(f)>2A\lambda) & \leq &
 \Phi\Big(\frac{A\lambda}{M}\Big)\sum_{k=j}^{N_A
 +j}\frac{1}{\Phi
 \Big(\frac{2^k}{M}\Big)}\\
 & + &\Phi\Big(\frac{A\lambda}{M}\Big)\sum_{k=N_A
 +j+1}^{\infty}\frac{1}{\Phi
 \Big(\frac{2^k}{M}\Big)}\\
&=&\textrm{I}+ \textrm{II}.
 \ee
 Now let's estimate \textrm{I} and
\textrm{II} respectively. Using Lemma 2.6 again, we get
 \be
\textrm{I} &= &\sum_{k=j}^{N_A
+j}\frac{\Phi\Big(\frac{A\lambda}{M}\Big)}{\Phi
\Big(\frac{2^k}{M}\Big)}\leq \sum_{k=j}^{N_A
+j}\Big(\frac{{A\lambda}}{{2^k}}\Big)^{\frac{1}{p}}
\\
& \leq & \sum_{k=j}^{N_A
+j}(2^{N_A +j+1-k})^{\frac{1}{p}} \leq \frac{2^{\frac{1}{p}(N_A +1)}}{1-2^{-\frac{1}{p}}}\\
& = & C_1  \ee
 and
 \be \textrm{II} &= &\sum_{k=N_A
 +j+1}^{\infty}\frac{\Phi\Big(\frac{A\lambda}{M}\Big)}{\Phi
 \Big(\frac{2^k}{M}\Big)}\leq \sum_{k=N_A
+j+1}^{\infty}\Big(\frac{{A\lambda}}{{2^k}}\Big)^{\frac{1}{q}}
 \\
& \leq & \sum_{k=N_A +j+1}^{\infty}\Big(2^{N_A +j+1-k}\Big)^{\frac{1}{q}}=\frac{1}{1-2^{-\frac{1}{q}}}\\
& = & C_2.  \ee
 Let $C_0=C_1+C_2$. It's easy to see that $C_0>1$. Thus
 \be\Phi\Bigg(\frac{2A\lambda}{(C_0\cdot
 c_\Phi)^{1/\ell}4 M}\Bigg)P\big(s(f)>2A\lambda\big) & \leq & c_\Phi
 \frac{1}{C_0\cdot c_\Phi}\Phi\Bigg(\frac{2A\lambda}{4
 M}\Bigg)P\big(s(f)>2A\lambda\big)\nonumber\\
 & \leq & c_\Phi \frac{1}{C_0\cdot c_\Phi}\cdot C_0 =1 .
 \ee
 And so we obtain
 \begin{equation}\|f\|_{wH_{\Phi}^{s}}\leq(C_0\cdot c_\Phi)^{1/\ell}4\cdot
 M.\end{equation}
 Consequently, (3.1) holds. The proof of Theorem 3.1 is complete.

 \begin{remark}
 If $f\in w\mathscr{H}_\Phi^s$ in Theorem
 3.1, then not only (i) and (ii) hold, but also
 the sum $\sum_{k=m}^n a^k$ converges to $f$ in $wH_\Phi^s$
as $m\rightarrow -\infty$, $n\rightarrow \infty$.
 Indeed,
 $$\sum_{k=m}^n a^k = \sum_{k=m}^n (f^{\nu_{k+1}}-f^{\nu_k})=f^{\nu_{n+1}}-f^{\nu_m}.$$
 By the sublinearity of $s(f)$ we have
 \be
 \|f-\sum_{k=m}^n a^k\|_{wH_\Phi^s}&=&\|s(f-f^{\nu_{n+1}}+f^{\nu_m})\|_{wL_\Phi}\leq
 \|s(f-f^{\nu_{n+1}})+s(f^{\nu_m})\|_{wL_\Phi}\\
 &\leq & (2c_\Phi)^{1/\ell}\cdot 2
 \bigg(\|s(f-f^{\nu_{n+1}})\|_{wL_\Phi}+\|s(f^{\nu_m})\|_{wL_\Phi}\bigg).
 \ee
 Since $s(f-f^{\nu_{n+1}})^2=s(f)^2- s(f^{\nu_{n+1}})^2$, then
 $s(f-f^{\nu_{n+1}})\leq s(f)$, $s(f^{\nu_m})\leq s(f)$
 and $s(f-f^{\nu_{n+1}}),\; s(f^{\nu_m})\rightarrow 0 $ a.e. as $m\rightarrow -\infty$, $n\rightarrow \infty$.
 Thus by Proposition 2.5, we have
 $$\|s(f-f^{\nu_{n+1}})\|_{wL_\Phi},\;\|s(f^{\nu_m})\|_{wL_\Phi} \rightarrow 0\;\;\;\; as\; m\rightarrow -\infty,\;n\rightarrow \infty,$$
 which means $\|f-\sum_{k=m}^n a^k\|_{wH_\Phi^s}\rightarrow 0$ as $m\rightarrow -\infty$, $n\rightarrow
 \infty$. Further, for $k\in \textbf{Z}$, $a^k=(a_n^k)_{n\geq 0}$ is $L_2$
 bounded, hence $H_2^s=L_2$ is dense in $w\mathscr{H}_\Phi^s.$
 \end{remark}

Recall that if $(\mathcal {F}_n)_{n\geq 0}$ is regular, then for any
 non-negative adapted sequence $\gamma=(\gamma_n)_{n\geq0}$ and $\lambda\geq
 \|\gamma_0\|_{\infty}$, there is a stopping time $\nu$ such that
 $$\{\gamma^*> \lambda\}\subset\{\nu < \infty\}, \quad \gamma_{\nu}^*\leq \lambda, \quad P(\nu<\infty)\leq R P(\gamma^* >\lambda)$$
 (see \cite{long}). Moreover, if $\lambda_1 \leq \lambda_2$, then we can take two stopping times $\nu_{\lambda_1}$ and $\nu_{\lambda_2}$ such
 that  $\nu_{\lambda_1}\leq \nu_{\lambda_2}$.
 Therefore, if $(\mathcal {F}_n)_{n\geq 0}$ is regular, we get the atomic decompositions for
 $wH_\Phi^S$ and $wH_\Phi$.

 \begin{thm}
 Let $\Phi \in \mathcal {G}_\ell$ with $\ell
\in (0,1]$ and $q_{\Phi^{-1}}<\infty$. Then, if $(\mathcal
{F}_n)_{n\geq 0}$ is regular, we have
 $$wH_\Phi^S= wH_{\Phi,at}^S\qquad with\;\,equivalent\;\,quasi\;\,norms;$$
 $$wH_\Phi= wH_{\Phi,at}\qquad with\;\,equivalent\;\,quasi\;\,norms.$$
 Moreover, if $f\in w\mathscr{H}_\Phi^S$ (or $w\mathscr{H}_\Phi$
 resp.), the sum  $\sum_{k=m}^n a^k$ converges to $f$ in $wH_\Phi^S$
 (or $wH_\Phi$ resp.) as $m\rightarrow -\infty$, $n\rightarrow
 \infty$.
 \end{thm}
 \noindent {\bf Proof.} The proof shall be given for $wH_\Phi^S$,
 only, since it is just slightly different from the one for
 $wH_\Phi$. Let $f\in wH_\Phi^S$. Then for sequence $S_n(f)$ and  $k\in
 \textbf{Z}$, take stopping times $\nu_k$ such that
 $$\{S(f)> 2^k\}\subset\{\nu_k < \infty\}, \quad S_{\nu_k}(f)\leq 2^k, \quad P(\nu_k<\infty)\leq R P(S(f) >2^k)$$
 and $\nu_k\leq \nu_{k+1}$, $\nu_k \uparrow \infty$. Still define $a_n^k=f_n^{\nu_{k+1}}-
 f_n^{\nu_{k}}$, then $a^k = (a_n^k)_{n\geq 0}$ is a martingale and
 \be
S(a^k)&=&\Big(\sum_{n=1}^{\infty}|da_n^k|^2\Big)^{1/2}= \Big(\sum_{n=1}^{\infty}|df_n\chi_{\{\nu_k<n\leq\nu_{k+1}\}}|^2\Big)^{1/2}\\
&\leq& S(f^{\nu_{k+1}})\leq 2^{k+1}=2\cdot2^k. \ee
 Thus,
 $(a_n^k)_{n\geq0}$ is $L_2$-bounded. Denote the limit still by
 $a^k$. Then $a_n^k=\E_n a^k$ for all $n\geq 0$. For $\nu_k\geq n$,
 $a_n^k = f_n^{\nu_{k+1}}- f_n^{\nu_k}= f_n-f_n=0.$ So $a^k$ is a
 w-2-atom and $f_n=\sum_{k\in \textbf{Z}}a_n^k$ a.e.
 Since $P(\nu_k<\infty)\leq R P(S(f) >2^k)$. Then let $C_0 =\max\{R,1\}$,
 we have
 \be
 \Phi\Bigg(\frac{2^k}{(c_\Phi C_0)^{1/\ell}\|f\|_{wH_{\Phi}^{S}}}\Bigg)P(\nu_k
 <\infty)&\leq&
 c_\Phi \cdot \frac{1}{c_\Phi C_0}\Phi\Bigg(\frac{2^k}{\|f\|_{wH_{\Phi}^{S}}}\Bigg)\cdot R P(s(f)>2^k)\\
  &\leq & \Phi\Bigg(\frac{2^k}{\|f\|_{wH_{\Phi}^{S}}}\Bigg)P(s(f)>2^k)\leq 1,\ee
 which means
 $$\inf \Big\{c>0:\sup_{k\in
 \textbf{Z}}\Phi\Big(\frac{2^k}{c}\Big)P(\nu_k <
 \infty)\leq 1\Big\}\leq (c_\Phi C_0)^{1/\ell} \|f\|_{wH_{\Phi}^{S}}<\infty.$$

 Conversely, suppose that $f\in wH_{\Phi,at}^S$. Let
$$M=\inf \Big\{c>0:\sup_{k\in \textbf{Z}}\Phi\Big(\frac{2^k}{c}\Big)P(\nu_k < \infty)\leq
1\Big\}.$$ Without loss of generality, here we also assume that
$A=2^{N_A}$, where $N_A\geq0$. For any $\lambda>0$, choose $j\in
\textbf{Z}$ such that $2^{j}\leq\lambda<2^{j+1}$. Define $f^N$ in
the same way as in (5.2), then similarly, we obtain
$$S(g)\leq \sum_{k=-\infty}^{j-1}S(a^k)\leq
A\lambda.$$
 And since $a_n^k=\E_n a^k = 0$ on the set $\{\nu_k\geq n\}$,
  $S(a^k)= 0$ on $\{\nu_k =\infty\}$. Thus,

\be \Phi\Big(\frac{2A\lambda}{4 M}\Big)P(S(f)>2A\lambda) & \leq &
\Phi\Big(\frac{A\lambda}{M}\Big)P(S(h)>A\lambda)\\
& \leq & \Phi\Big(\frac{A\lambda}{M}\Big)\sum_{k=j}^{\infty}P(\nu_k
<\infty) \ee
 Dealing with the last inequality in the same way as in Theorem 3.1, we
 obtain
 $$\|f\|_{wH_{\Phi}^{S}}\leq C
 M.$$

 Further, if $f\in w\mathscr{H}_\Phi^S$, just as Remark 3.2, by Proposition 2.5 we
 get that the sum  $\sum_{k=m}^n a^k$ converges to $f$ in $wH_\Phi^S$
 as $m\rightarrow -\infty$, $n\rightarrow\infty.$
The proof of Theorem 3.2 is complete.

\begin{thm} Let $\Phi \in \mathcal {G}_\ell$ with $\ell
\in (0,1]$ and $q_{\Phi^{-1}}<\infty$. Then
 $$w\mathcal {Q}_\Phi= wH_{\Phi,at}^S\qquad with\;\,equivalent\;\,quasi\;\,norms;$$
 $$w\mathcal {D}_\Phi = wH_{\Phi,at}\qquad with\;\,equivalent\;\,quasi\;\,norms.$$
\end{thm}

 \noindent {\bf Proof.} The proof of Theorem 3.4 is similar to
 that of Theorem 3.1. So we just sketch the outline and omit the
 details. Suppose that $f =(f_n)_{n\geq 0} \in
 w\mathcal {Q}_\Phi$ (or $f =(f_n)_{n\geq 0} \in
 w\mathcal {D}_\Phi$ resp.). Let $\nu_k = \inf\{n:
 \lambda_n>2^k\}$, where
 $(\lambda_n)_{n\geq0}$ is the sequence in the definition of $w\mathcal {Q}_\Phi$
 (or $w\mathcal {D}_\Phi$, resp.) For $k\in\textbf{Z}$
 and $n\geq 0 $, we still define $a_n^k=f_n^{\nu_{k+1}}-
 f_n^{\nu_{k}}$. Then,  in the same way as in Theorem
 3.1, we can prove that there exists $A>0$ such that $S(a^k)\leq
 A \cdot 2^k$ (or $M(a^k)\leq A\cdot
 2^k$ resp.), and that $\|f\|_{wH_{\Phi,at}^S}\leq C
 \|f\|_{w\mathcal {Q}_\Phi}$ (or $\|f\|_{wH_{\Phi,at}}\leq C
 \|f\|_{w\mathcal {D}_\Phi}$, resp.).

 For the converse part, assume that $f =(f_n)_{n\geq 0} \in wH_{\Phi,at}^{S}$ (or $f =(f_n)_{n\geq 0} \in wH_{\Phi,at}$
 resp.), and let $\lambda_n=\sum_{k\in \textbf{Z}}\chi(\nu_k \leq
 n)\|S(a^k)\|_\infty$ (or $\lambda_n=\sum_{k\in \textbf{Z}}\chi(\nu_k
 \leq n)\|M(a^k)\|_\infty$ resp.). Then $(\lambda_n)_{n\geq 0}$ is a
 non-negative, non-decreasing and adapted sequence with
 $S_{n+1}(f)\leq \lambda_n$ (or $M_{n+1}(f)\leq \lambda_n$ resp.).
 For $y>0$, choose $j\in \textbf{Z}$ such that $2^j\leq y
 <2^{j+1}$. Then $\lambda_\infty =
 \lambda_\infty^{(1)}+\lambda_\infty^{(2)}$ with
 $\lambda_\infty^{(1)}=\sum_{k=-\infty}^{j-1} \chi(\nu_k \leq
 n)\|S(a^k)\|_\infty$ and  $\lambda_\infty^{(2)}=\sum_{k=j}^{\infty}\chi(\nu_k \leq
 n)\|S(a^k)\|_\infty$ (or, $\lambda_\infty^{(1)}=\sum_{k=-\infty}^{j-1} \chi(\nu_k \leq
 n)\|M(a^k)\|_\infty$ and  $\lambda_\infty^{(2)}=\sum_{k=j}^{\infty}\chi(\nu_k \leq
 n)\|M(a^k)\|_\infty$ resp.).
 Similar to the argument of (3.3) (replacing $s(g)$ and $s(h)$ by $\lambda_\infty^{(1)}$
 and $\lambda_\infty^{(2)}$,
 resp.), we obtain $\|f\|_{w\mathcal {Q}_\Phi}\leq C
 \|f\|_{wH_{\Phi,at}^S}$ (or $\|f\|_{w\mathcal {D}_\Phi}\leq C
 \|f\|_{wH_{\Phi,at}}$ resp.).\\

 \noindent Theorem 3.3 and Theorem 3.4 together give the following corollary.
\begin{cor}
 Let $\Phi \in \mathcal {G}_\ell$ with $\ell
 \in (0,1]$ and $q_{\Phi^{-1}}<\infty$. If $(\mathcal {F}_n)_{n\geq
 0}$ is regular, then
 $$w\mathcal {Q}_\Phi= wH_\Phi^S,\qquad w\mathcal {D}_\Phi= wH_\Phi.$$
 \end{cor}

%%%%%%%%%%%%%%%%%%%%%%%%%%%%%%%%%%%%%%%%%%%%%%%%%%%%%%%%%%%%%%%%
%%%%%%%%%%%%%%%%%%%%%%%%%%%%%%%%%%%%%%%%%%%%%%%%%%%%%%%%%%%%%%%%

\section{Bounded operators on weak Orlicz-Hardy spaces}

%%%%%%%%%%%%%%%%%%%%%%%%%%%%%%%%%%%%%%%%%%%%%%%%%%%%%%%%%%%%%%%%
%%%%%%%%%%%%%%%%%%%%%%%%%%%%%%%%%%%%%%%%%%%%%%%%%%%%%%%%%%%%%%%%

\smallskip

As one of the applications of the atomic decompositions, we shall
obtain a sufficient condition for a sublinear operator to be bounded
from the weak martingale Orlicz-Hardy spaces to $wL_\Phi$ spaces.
Applying the condition to $M(f)$, $S(f)$ and $s(f)$, we deduce a
series of martingale inequalities.

An operator $T: X \rightarrow Y$ is called a sublinear operator if
it satisfies
$$|T(f+g)|\leq |Tf|+|Tg|,\;\; |T(\alpha f)|\leq |\alpha||Tf|,$$
where $X$ is a martingale spaces, $Y$ is a measurable function
space.

\begin{thm} Let $1\leq r\leq 2$ and $T: L_r(\Omega) \rightarrow L_r(\Omega)$
be a bounded sublinear operator. If
\begin{equation}P(|Ta|>0)\leq C P(\nu<\infty)\end{equation}
for all w-1-atoms, where $\nu$ is the stopping time associated with
$a$ and $C$ is a positive constant, then, for $\Phi\in \mathcal
{G}_\ell$ with $q_{\Phi^{-1}}<\infty$ and $1/p_{\Phi^{-1}}<r$, there
exists a positive constant $C$ such that
$$\|Tf\|_{wL_\Phi}\leq C\|f\|_{wH_\Phi^s},\quad f\in wH_\Phi^s.$$

\end{thm}

\noindent {\bf Proof.} Assume that $f\in wH_\Phi^s$. By Theorem 3.1,
$f$ can be decomposed into the sum of a sequence of w-1-atoms. For
any fixed $\lambda > 0$, choose $j\in \textbf{Z}$ such that
$2^{j-1}\leq \lambda < 2^j$ and let
$$f=\sum_{k=-\infty}^{\infty}a^k=\sum_{k=-\infty}^{j-1}a^k+\sum_{k=j}^{\infty}a^k:=g+h.$$
It follows from the sublinearity of $T$ that $|Tf|\leq|Tg|+|Th|$, so
$$P(|Tf|>2\lambda) \leq
P(|Tg|>\lambda)+P(|Th|>\lambda).$$ In Theorem 3.1, we have proved
that $s(a^k)\leq A\cdot 2^k$ for some $A>0$ and $s(a^k)=0$ on the
set $\{\nu_k=\infty\}$. Denote $p=p_{\Phi^{-1}}$, $q=q_{\Phi^{-1}}$.
Remember that
\begin{equation}\label{3}\|a^k\|_r\leq C\|s(a^k)\|_r,\quad \quad 1\leq r\leq 2.
\end{equation}  It results from Lemma 2.6 that
\be \|g\|_r & \leq & \sum_{k=-\infty}^{j-1}\|a^k\|_r\leq
C\sum_{k=-\infty}^{j-1}\|s(a^k)\|_r\\
&\leq & C\sum_{k=-\infty}^{j-1} 2^k P(\nu_k <\infty)^{\frac{1}{r}}.\\
\ee Since $T$ is bounded on $L_r(\Omega)$, then
 \be
\Phi\Big(\frac{\lambda}{\|f\|_{wH_\Phi^s}}\Big)P(|T g|>\lambda) &
\leq &\Phi\Big(\frac{\lambda}{\|f\|_{wH_\Phi^s}}\Big)
\frac{\|T g\|_r^r}{\lambda^r}\\
& \leq & C \Phi\Big(\frac{\lambda}{\|f\|_{wH_\Phi^s}}\Big)
\frac{\|g\|_r^r}{\lambda^r}.\ee By the estimate of $\|g\|_r$ above,
we get
 \be \Phi\Big(\frac{\lambda}{\|f\|_{wH_\Phi^s}}\Big)P(|T
g|>\lambda) & \leq & C
\Phi\Big(\frac{\lambda}{\|f\|_{wH_\Phi^s}}\Big)\Bigg(\sum_{k=-\infty}^{j-1}\frac{2^kP(\nu_k
<\infty)^{\frac{1}{r}}}{\lambda}\Bigg)^r \\
& = &C
\Phi\Big(\frac{\lambda}{\|f\|_{wH_\Phi^s}}\Big)\Bigg(\sum_{k=-\infty}^{j-1}\frac{2^k
\Phi\Big(\frac{2^k}{\|f\|_{wH_\Phi^s}}\Big)^{\frac{1}{r}}P(\nu_k
<\infty)^{\frac{1}{r}}}{\lambda \Phi\Big(\frac{2^k}{\|f\|_{wH_\Phi^s}}\Big)^{\frac{1}{r}}}\Bigg)^r\\
& \leq & C \Bigg(\sum_{k=-\infty}^{j-1}\frac{2^k
\Phi\Big(\frac{\lambda}{\|f\|_{wH_\Phi^s}}\Big)^{\frac{1}{r}}}{\lambda\Phi\Big(\frac{2^k}{\|f\|_{wH_\Phi^s}}\Big)^{\frac{1}{r}}}\Bigg)^r.\ee
Using Lemma 2.6 and the condition
$\frac{1}{p}=\frac{1}{{p_{\Phi^{-1}}}}<r$, we obtain \be
\Phi\Big(\frac{\lambda}{\|f\|_{wH_\Phi^s}}\Big)P(|T g|>\lambda) &
\leq & C \Bigg(\sum_{k=-\infty}^{j-1}\Big(\frac{2^k}{\lambda}\Big)
\Big(\frac{\lambda}{2^k}\Big)^{\frac{1}{pr}}\Bigg)^r\\
& = & C \lambda^{\frac{1}{p}-r}\cdot
\Bigg(\sum_{k=1-j}^{\infty}\Big(\Big(\frac{1}{2}\Big)^{1-\frac{1}{pr}}\Big)^k
\Bigg)^r \\
&\leq & C \lambda^{\frac{1}{p}-r}\cdot (2^{j-1})^{r
 -\frac{1}{p}}\leq C_1.\ee
 Taking $C_{\textrm{I}}=\big(2c_\Phi\max\{C_1, 1\}\big)^{1/\ell}$,
 then
 \be\Phi\Big(\frac{\lambda}{C_\textrm{I}\|f\|_{wH_\Phi^s}}\Big)P(|T
g|>\lambda) & \leq & c_\Phi \frac{1}{2c_\Phi\max\{C_1,
1\}}\Phi\Big(\frac{\lambda}{\|f\|_{wH_\Phi^s}}\Big)P(|T
g|>\lambda)\\
 & \leq &
  c_\Phi \frac{1}{2c_\Phi\max\{C_1,
 1\}} C_1 \leq \frac{1}{2}.\ee
 On the other hand, since $|Th|\leq \sum_{k=j}^{\infty}|Ta^k|$, we
 get
 \be
 \Phi\Big(\frac{\lambda}{\|f\|_{wH_\Phi^s}}\Big)P(|T
 h|>\lambda) & \leq &
 \Phi\Big(\frac{\lambda}{\|f\|_{wH_\Phi^s}}\Big)P(|T
 h|>0)\\
 & \leq &
 \Phi\Big(\frac{\lambda}{\|f\|_{wH_\Phi^s}}\Big)\sum_{k=j}^{\infty}P(|T
 a^k|>0)\\
 & \leq &
 C \Phi\Big(\frac{\lambda}{\|f\|_{wH_\Phi^s}}\Big)\sum_{k=j}^{\infty}P(\nu_k<\infty)\\
 & = & C \sum_{k=j}^{\infty}\frac{\Phi\Big(\frac{\lambda}{\|f\|_{wH_\Phi^s}}\Big)}{\Phi\Big(\frac{2^k}{\|f\|_{wH_\Phi^s}}\Big)}\Phi\Big(\frac{2^k}{\|f\|_{wH_\Phi^s}}\Big)P(\nu_k<\infty)\\
 &\leq & C
 \sum_{k=j}^{\infty}\frac{\Phi\Big(\frac{\lambda}{\|f\|_{wH_\Phi^s}}\Big)}{\Phi\Big(\frac{2^k}{\|f\|_{wH_\Phi^s}}\Big)}\leq
 C
 \sum_{k=j}^{\infty}\Big(\frac{\lambda}{2^k}\Big)^{\frac{1}{q}}\leq
 C_2. \ee
 Taking $C_{\textrm{II}}=\big(2c_\Phi\max\{C_2,
 1\}\big)^{1/\ell}$,
 then
 $$\Phi\Big(\frac{\lambda}{C_\textrm{II}\|f\|_{wH_\Phi^s}}\Big)P(\Phi(|T
 h|)>\lambda) \leq \frac{1}{2}.$$
 Since $T$ is sublinear,
 \be
 \Phi\Big(\frac{2\lambda}{2(C_\textrm{I}+C_\textrm{II})\|f\|_{wH_\Phi^s}}\Big)P(|T
 f|>2\lambda) &\leq &
 \Phi\Big(\frac{\lambda}{C_\textrm{I}\|f\|_{wH_\Phi^s}}\Big)P(|T
 g|>\lambda) \\
 &+&
 \Phi\Big(\frac{\lambda}{C_\textrm{II}\|f\|_{wH_\Phi^s}}\Big)P(|T
 h|>\lambda) \\
 &\leq & 1 .\ee
 Hence,
 $$\|Tf\|_{wL_\Phi}\leq C\|f\|_{wH_\Phi^s},\quad f\in wH_\Phi^s.$$
 The proof is complete.

\begin{remark}
Similarly, if $T: L_r(\Omega) \rightarrow L_r(\Omega)$ is a bounded
sublinear operator, $1\leq r<\infty$, and (4.1) holds for all
w-2-atoms (or w-3-atoms). Then for
 $\Phi\in \mathcal {G}_\ell$ with
$q_{\Phi^{-1}}<\infty$ and $\frac{1}{{p_{\Phi^{-1}}}}< r$, there
exists a constant $C>0$ such that
$$\|Tf\|_{wL_\Phi}\leq C\|f\|_{w\mathcal {Q}_\Phi},\quad f\in w\mathcal {Q}_\Phi$$
$$(or\;\; \|Tf\|_{wL_\Phi}\leq C\|f\|_{w\mathcal {D}_\Phi},\quad f\in w\mathcal {D}_\Phi).$$
\noindent In this case, we do not need to restrict $1\leq r\leq2$;
in fact, \eqref{3} is replaced by $\|a^k\|_r\leq C\|S(a^k)\|_r$(or
$\|a^k\|_r\leq C\|M(a^k)\|_r$), which always holds for $1\leq
r<\infty$ by the Burkholder-Davis-Gundy inequalities.
\end{remark}

\begin{thm}
 Let $\Phi\in \mathcal {G}_\ell$ with $\ell\in(0,1]$ and $q_{\Phi^{-1}}<\infty$. Then for all
 martingales $f=(f_n)_{n\geq 0}$ the following inequalities hold:
\begin{equation}\label{4}\|f\|_{wH_\Phi}\leq C\|f\|_{wH_\Phi^s},\quad  \|f\|_{wH_\Phi^S}\leq C\|f\|_{wH_\Phi^s};\end{equation}
\begin{equation}\label{5}\|f\|_{wH_\Phi}\leq C\|f\|_{w\mathcal {Q}_\Phi},\quad  \|f\|_{wH_\Phi^S}\leq C\|f\|_{w\mathcal {Q}_\Phi},\quad \|f\|_{wH_\Phi^s}\leq C\|f\|_{w\mathcal {Q}_\Phi};\end{equation}
\begin{equation}\label{6}\|f\|_{wH_\Phi}\leq C\|f\|_{w\mathcal {D}_\Phi},\quad  \|f\|_{wH_\Phi^S}\leq C\|f\|_{w\mathcal {D}_\Phi},\quad \|f\|_{wH_\Phi^s}\leq C\|f\|_{w\mathcal {D}_\Phi};\end{equation}
\begin{equation}\label{7}C^{-1}\|f\|_{w\mathcal {D}_\Phi} \leq \|f\|_{w\mathcal {Q}_\Phi} \leq C\|f\|_{w\mathcal {D}_\Phi}.\end{equation}
Moreover, if $\{\mathcal {F}_n\}_{n\geq 0}$ is regular, then
$wH_\Phi^S=w\mathcal {Q}_\Phi=w\mathcal {D}_\Phi=wH_\Phi =
wH_\Phi^s.$

\end{thm}

\noindent {\bf Proof.} First we show \eqref{4}. Let $f\in
wH_\Phi^s$. The maximal operator $Tf=Mf$ is sublinear. It's well
know that $T$ is $L_2$-bounded. If $a$ is a w-1-atom and $\nu$ is
the stopping time associated with $a$, then
$\{|Ta|>0\}=\{M(a)>0\}\subset\{\nu<\infty\}$ and hence (4.1) holds.
Since $\Phi\in \mathcal {G}_\ell $, the condition
$\frac{1}{{p_{\Phi^{-1}}}}< 2$ always holds for convex function
$\Phi^{-1}$. Thus it follows from Theorem 4.1 that
$$\|f\|_{wH_\Phi}= \|Tf\|_{wL_\Phi}\leq C\|f\|_{wH_\Phi^s}.$$
Similarly, considering the operator $Tf=Sf$ we get the second
inequality of \eqref{4}.

Next we show \eqref{5} and \eqref{6}. Choose $r$ such that $
\frac{1}{p_{\Phi^{-1}}}<2 < r< \infty$. Noticing that the operator
$Mf$, $Sf$ and $sf$ are $L_r$ bounded. Taking $Tf=Mf$, $Sf$ or $sf$,
resp. By Remark 4.2, we get \eqref{5} and \eqref{6}.

To prove \eqref{7}, we use \eqref{5} and \eqref{6}. The method used
below is the same as the proof of Theorem 3.5 in \cite{ren2}. Assume
that $f=(f_n)_{n\geq 0}\in w\mathcal {Q}_\Phi$, then there exists an
optimal control $(\lambda_n^{(1)})_{n\geq 0}$ such that $S_n(f)\leq
\lambda_{n-1}^{(1)}$ with $\lambda_\infty^{(1)}\in wL_\Phi$. Since
$$|f_n|\leq f_{n-1}^* +\lambda_{n-1}^{(1)},$$ by \eqref{5} we have
$$\|f\|_{w\mathcal {D}_\Phi}\leq C (\|f\|_{wH_\Phi}+\|\lambda_\infty^{(1)}\|_{wL_\Phi})\leq C\|f\|_{w\mathcal {Q}_\Phi}.$$
On the other hand, if $f=(f_n)_{n\geq 0}\in w\mathcal {D}_\Phi$,
then there exists an optimal control $(\lambda_n^{(2)})_{n\geq 0}$
such that $|f_n|\leq \lambda_{n-1}^{(2)}$ with
$\lambda_\infty^{(2)}\in wL_\Phi$. Notice that $$S_n(f)\leq
S_{n-1}(f)+2\lambda_{n-1}^{(2)}.$$ Using \eqref{6} we can get the
other side of  \eqref{7}.

Further, suppose that $\{\mathcal {F}_n\}_{n\geq 0}$ is regular.
Then for any martingale $f=(f_n)_{n\geq 0}$, we have $|df_n|^2 \leq
\frac{R-1}{2} \E_{n-1}|df_n|^2$ (see \cite{weisz3}, pp 31,
Proposition 2.19). Thus
$$S_n(f)\leq \sqrt{\frac{R-1}{2}}s_n(f).$$
Since $s_n(f)\in \mathcal {F}_{n-1}$, by the definition of
$w\mathcal {Q}_\Phi$ we have
$$\|f\|_{w\mathcal {Q}_\Phi}\leq \|s(f)\|_{wL_\Phi}=\|f\|_{wH_\Phi^s}.$$
Using \eqref{5}, \eqref{7} and Corollary 3.5, we conclude that
$$wH_\Phi^S=w\mathcal {Q}_\Phi=w\mathcal {D}_\Phi=wH_\Phi =
wH_\Phi^s.$$

%%%%%%%%%%%%%%%%%%%%%%%%%%%%%%%%%%%%%%%%%%%%%%%%%%%%%%%%%%%%%%%%
%%%%%%%%%%%%%%%%%%%%%%%%%%%%%%%%%%%%%%%%%%%%%%%%%%%%%%%%%%%%%%%%

\section{The duality results}

%%%%%%%%%%%%%%%%%%%%%%%%%%%%%%%%%%%%%%%%%%%%%%%%%%%%%%%%%%%%%%%%
%%%%%%%%%%%%%%%%%%%%%%%%%%%%%%%%%%%%%%%%%%%%%%%%%%%%%%%%%%%%%%%%

\smallskip

In this section, we investigate the dual of weak martingale
Orlicz-Hardy spaces and give a new John-Nirenberg theorem.

\begin{thm} Let $\Phi \in \mathcal {G}_\ell$ with $\ell
\in (0,1]$, $q_{\Phi^{-1}}<\infty$ and
$\phi(r)=1/(r\Phi^{-1}(1/r))$. Then
$$(w\mathscr{H}_\Phi^s)^*=w\mathcal {L}_{2,\phi}.$$
\end{thm}

\noindent {\bf Proof.} Let $g\in w\mathcal {L}_{2,\phi}$, then $g\in
H_2^s$. Define $$l_g(f)=\E\Big(\sum_{n=1}^{\infty}df_n
dg_n\Big),\qquad f\in H_2^s.$$ From Theorem 3.1, there is a sequence
of w-1-atoms $a^k$ and corresponding stopping times $\nu_k$, where
$k\in \textbf{Z}$, such that
$$s(a^k)\leq 2^{k+1},\quad \Phi\Bigg(\frac{2^k}{\|f\|_{wH_{\Phi}^{s}}}\Bigg)P(\nu_k
<\infty)\leq 1$$ and
$$df_n =\sum_{k=-\infty}^{\infty}da_n^k \quad \textrm{a.e.}$$
 for all $n\in \textbf{N}$. The last series
also converges to $df_n$ in $H_2^s$-norm. Hence
$$l_g(f)=\sum_{n=1}^{\infty}\sum_{k=-\infty}^{\infty}\E (da_n^k dg_n).$$
Applying the H\"{o}lder inequality and the definition of weak atoms,
we get that \be |l_g(f)| &\leq & \sum_{k=-\infty}^{\infty}\E
\Big(\sum_{n=1}^{\infty}|da_n^k|\chi_{\{\nu_k<n\}}|dg_n|\Big)\\
&\leq& \sum_{k=-\infty}^{\infty} \Big(\E
\sum_{n=1}^{\infty}|da_n^k|^2\Big)^{1/2}\Big(\E
\sum_{n=1}^{\infty}|dg_n|^2 \chi_{\{\nu_k<n\}}\Big)^{1/2}\\
&=& \sum_{k=-\infty}^{\infty} \Big(\E
\sum_{n=1}^{\infty}|da_n^k|^2\Big)^{1/2}\|S(g-g^{\nu_k})\|_2. \ee
Since $\Big(\E \sum_{n=1}^{\infty}|da_n^k|^2\Big)^{1/2}=\Big(\E (S^2
(a^k))\Big)^{1/2}=\|s(a^k)\|_2\leq 2^{k+1}P(\nu_k <\infty)^{1/2}$
and $P(\nu_k<\infty)\leq
1\Big/\Phi\Big(\frac{2^k}{\|f\|_{wH_{\Phi}^{s}}}\Big)$, then
 \be
 |l_g(f)| &\leq & \sum_{k=-\infty}^{\infty}2^{k+1}P(\nu_k
 <\infty)^{1/2}\|g-g^{\nu_k}\|_2\\
 &\leq&\sum_{k=-\infty}^{\infty}2^{k+1}\Bigg(\Phi\Big(\frac{2^k}{\|f\|_{wH_{\Phi}^{s}}}\Big)\Bigg)^{-1/2}\|g-g^{\nu_k}\|_2.
 \ee Let
 $A_k=1\Big/\Phi\Big(\frac{2^k}{\|f\|_{wH_{\Phi}^{s}}}\Big)$
and still denote $p=p_{\Phi^{-1}}$, $q=q_{\Phi^{-1}}$. By Lemma 2.6,
we obtain
 \be |l_g(f)| &\leq &
2\|f\|_{wH_{\Phi}^{s}}\sum_{k=-\infty}^{\infty}\frac{2^k}{\|f\|_{wH_{\Phi}^{s}}}\cdot
A_k ^{1/2}\|g-g^{\nu_k}\|_2\\
&\leq & 2
\|f\|_{wH_{\Phi}^{s}}\sum_{k=-\infty}^{\infty}\frac{1}{\phi(A_k)}A_k
^{-1/2}\sup_{P(\nu_k <\infty)\leq A_k}\|g-g^{\nu_k}\|_2\\
&=& 2^q
\|f\|_{wH_{\Phi}^{s}}\sum_{k=-\infty}^{\infty}t_{\phi}^{2}(A_k). \ee
Using Lemma 2.6 again, we get $$ \frac{A_{k+1}}{A_k}=
\frac{\Phi\Big(\frac{2^k}{\|f\|_{wH_\Phi^s}}\Big)}{\Phi\Big(\frac{2^{k+1}}{\|f\|_{wH_\Phi^s}}\Big)}
\leq \Big(\frac{2^k}{2^{k+1}}\Big)^{1/q}=
\Big(\frac{1}{2}\Big)^{1/q}.$$ Thus, \be
\sum_{k=-\infty}^{\infty}t_{\phi}^{2}(A_k) & = &
\sum_{k=-\infty}^{\infty}\frac{t_{\phi}^{2}(A_k)(A_k - A_{k+1})}{A_k
- A_{k+1}}\leq
\frac{1}{1-(\frac{1}{2})^{1/q}}\sum_{k=-\infty}^{\infty}\frac{t_{\phi}^{2}(A_k)(A_k
- A_{k+1})}{A_k }\\
&\leq & C\int_{0}^{\infty}\frac{t_{\phi}^{2}(x)}{x}dx
=C\|g\|_{w\mathcal {L}_{2,\phi}}. \ee And so $$|l_g(f)| \leq
C\|f\|_{wH_{\Phi}^{s}}\|g\|_{w\mathcal {L}_{2,\phi}}.$$
 Since $H_2^s$ is dense in $w\mathscr{H}_\Phi^s$ (see Remark 3.2), $l$ can be uniquely extended to a
 continuous linear functional on $w\mathscr{H}_\Phi^s$.

 Conversely, let $l\in (w\mathscr{H}_\Phi^s)^*$. Note that $H_2^s \subset
 w\mathscr{H}_\Phi^s$, hence $l \in (H_{2}^{s})^*$. That means there exists
 $g\in H_2^s$ such that
 $$l(f)=\E\Big(\sum_{n=1}^{\infty}df_n dg_n\Big),\qquad f\in H_2^s.$$
 Let $\nu_k$ be the stopping times satisfying $P(\nu_k <\infty)\leq
 2^{-k}$ ($k\in \textbf{Z}$). For $k\in \textbf{Z}$, we define
 $$a^k=\frac{g-g^{\nu_k}}{(2^{k})^{1/2}\frac{1}{\Phi^{-1}(2^k)}\|s(g-g^{\nu_k})\|_2}.$$
 The function $a^k$ is not necessarily a weak atom. However, it
 satisfies the condition (i) of Definition 2.8, namely, $a_n^k=0$ on
 the set $\{\nu_k\geq n\}$. For $\lambda>0$, choose $j\in \textbf{Z}$
 such that $2^{j}\leq \lambda < 2^{j+1}$ and define the martingales
 $f^N$, $g^N$ and $h^N$, respectively, by
 \begin{equation}f_n^N=\sum_{k=-N}^{N}a_n^k,\quad g_n^N=\sum_{k=-N}^{j-1}a_n^k \quad \textrm{and} \quad
 h_n^N=\sum_{k=j}^{N}a_n^k.\end{equation}
 Then we have $\Phi(s(f))\leq\Phi(s(g))+\Phi(s(h))$ by the
sublinearity of $s(f)$ and $\Phi(t)$. And thus
$$P(\Phi(s(f^N))>2\lambda)\leq P(\Phi(s(g^N))>\lambda)+P(\Phi(s(h^N))>\lambda).$$
 Since $$\|s(g^N)\|_2 \leq \sum_{k=-N}^{j-1}\|s(a^k)\|_2\leq
 \sum_{k=-N}^{j-1}(2^{-k})^{1/2}\Phi^{-1}(2^k) ,$$ then \be
 P(s(g^N)>\Phi^{-1}(\lambda)) &\leq
 &\frac{1}{(\Phi^{-1}(\lambda))^2}\|s(g^N)\|_2^2\\
 &\leq & \frac{1}{(\Phi^{-1}(\lambda))^2}
 \Bigg(\sum_{k=-N}^{j-1}(2^{-k})^{1/2}\Phi^{-1}(2^k)\Bigg)^{2}\\
 &=&\Bigg(\sum_{k=-N}^{j-1}\frac{(2^{-k})^{1/2}\Phi^{-1}(2^k)}{\Phi^{-1}(\lambda)}\Bigg)^2\\
 &\leq &\Bigg(\sum_{k=-N}^{j-1}(2^{-k})^{1/2}\Big(\frac{2^k}{\lambda}\Big)^p\Bigg)^2\\
 &=&\lambda^{-2p}\Big(\sum_{k=-N}^{j-1}(2^{p-\frac{1}{2}})^k\Big)^2\leq
 C_1\lambda^{-1}. \ee
 In the last inequality above, we used $1/2<1\leq p \leq q <\infty$, which  results from that $\Phi^{-1}$
 is a convex function.
 Denote
 $C_\textrm{I}=(2c_{\Phi}\max\{C_1,1\})^{1/\ell}$, then
 $$\Phi\Big(\frac{\Phi^{-1}(\lambda)}{C_\textrm{I}}\Big)P(s(g^N)>\Phi^{-1}(\lambda))\leq c_\Phi\frac{1}{2c_{\Phi}\max\{C_1,1\}}\lambda P(s(g^N)>\Phi^{-1}(\lambda))\leq \frac{1}{2}.$$
 Noticing that $a_n^k=0$ on $\{\nu_k\geq n\}$ and $P(\nu_k
 <\infty)\leq 2^{-k}$, we get
 $$P(s(h^N)>\Phi^{-1}(\lambda))\leq \sum_{k=j}^{N}P(\nu_k<\infty)\leq \sum_{k=j}^{N}2^{-k}=2^{1-j}\leq 4\lambda^{-1}.$$
 Denote $C_{\textrm{II}}=(8c_{\Phi})^{1/\ell}$, then
 $$\Phi\Big(\frac{\Phi^{-1}(\lambda)}{C_\textrm{II}}\Big)P(s(h^N)>\Phi^{-1}(\lambda))\leq \frac{1}{2}.$$
 Let $C=2^q\max\{C_\textrm{I}, C_{\textrm{II}}\}$, then
 \be
 \Phi\Big(\frac{\Phi^{-1}(2\lambda)}{C}\Big)P(s(f^N)>\Phi^{-1}(2\lambda))&\leq&
 \Phi\Big(\frac{\Phi^{-1}(2\lambda)}{C}\Big)\Big(P(s(g^N)>\Phi^{-1}(\lambda))+P(s(h^N)>\Phi^{-1}(\lambda))\Big)\\
 &\leq &\Phi\Big(\frac{2^q\Phi^{-1}(\lambda)}{2^q
 C_\textrm{I}}\Big)P(s(g^N)>\Phi^{-1}(\lambda))\\
 &+& \Phi\Big(\frac{2^q\Phi^{-1}(\lambda)}{2^q
 C_\textrm{II}}\Big)P(s(h^N)>\Phi^{-1}(\lambda))\\
 &\leq& \frac{1}{2}+ \frac{1}{2}=1, \ee which implies $\|f^N\|_{
 wH_{\Phi}^{s}} \leq C$. Since \be
 l(f^N)&= &\E\sum_{n=1}^{\infty}df_n^N dg_n=\E\sum_{n=1}^{\infty}\sum_{k=-N}^{N}da_n^k dg_n\\
 &=& \sum_{k=-N}^{N}\E
 \sum_{n=1}^{\infty}\frac{|dg_n-dg_n^{\nu_k}|^2}{(2^k)^{1/2}\Phi^{-1}(2^k)\|S(g-g^{\nu_k})\|_2}\\
 &=& \sum_{k=-N}^{N}
(2^{-k})^{1/2}\frac{1}{\Phi^{-1}(2^k)}\|g-g^{\nu_k}\|_2, \ee then
 $$C\|l\|\geq l(f^N)=\sum_{k=-N}^{N}
 \frac{1}{\phi(2^{-k})}(2^{-k})^{-1/2}\|g-g^{\nu_k}\|_2.$$ Taking
 over all $N\in \textbf{N}$ and the supremum over all of such
 stopping times such that $P(\nu_k <\infty)\leq 2^{-k}$, $k\in
 \textbf{Z}$, we obtain
 $$\|g\|_{w\mathcal {L}_{2,\phi}}=\int_{0}^{\infty}\frac{t_{\phi}^{2}(x)}{x}dx \leq C\sum_{k=-\infty}^{\infty}t_{\phi}^{2}(2^{-k})\leq C\|l\|.$$
 The proof of Theorem 5.1 is complete.\\

 To obtain the new John-Nirenberg theorem, we first prove two lemmas.
 \begin{lem} Let $\Phi \in \mathcal {G}_\ell$ with $\ell
 \in (0,1]$, $q_{\Phi^{-1}}<\infty$ and
 $\phi(r)=1/(r\Phi^{-1}(1/r))$. If $\{\mathcal {F}_n\}_{n\geq 0}$ is
 regular, then
 $$(w\mathscr{H}_\Phi)^*=w\mathcal {L}_{1,\phi}.$$

 \end{lem}

 \noindent {\bf Proof.} Let $g\in w\mathcal {L}_{1,\phi}$ and define
 $$l_g(f)= \E(fg),\qquad f\in L_\infty.$$ Then
 \be
 |l_g(f)|& = &|\E(f g)|= |\sum_{k=-\infty}^\infty
 \E(a^k(g-g^{\nu_k}))|\\
 & \leq & \sum_{k=-\infty}^\infty \|a^k\|_\infty
 \|g-g^{\nu_k}\|_1 \leq \sum_{k=-\infty}^\infty 2^{k+1}
 \|g-g^{\nu_k}\|_1\\
 & \leq & 2 \|f\|_{wH_\Phi}\sum_{k=-\infty}^\infty \frac{ 2^k}{\|f\|_{wH_\Phi}}\|g-g^{\nu_k}\|_1.
 \ee
 Denote $A_k=1\Big/\Phi\Big(\frac{2^k}{\|f\|_{wH_\Phi}}\Big)$, we get
 \be
 |l_g(f)| & \leq &  2\|f\|_{wH_\Phi}\sum_{k=-\infty}^\infty
 \frac{1}{\phi(A_k)} A_k^{-1}\|g-g^{\nu_k}\|_1 \\
 & \leq & 2 \|f\|_{wH_\Phi}\sum_{k=-\infty}^\infty t_\phi^1(A_k)\\
 & \leq & C\|f\|_{wH_\Phi}\|g\|_{w\mathcal {L}_{1,\phi}}
 \ee

 Conversely, suppose that $l\in (w\mathscr{H}_\Phi)^* $. Since $L_2$
 is dense in $w\mathscr{H}_\Phi$,
 there exists $g\in L_2 \subset L_1$ such that
 $$l(f)= \E(fg),\qquad f\in L_\infty.$$
 Let $\nu_k$ be the stopping times satisfying $P(\nu_k <\infty)\leq
 (\Phi(2^k))^{-1}$ ($k\in \textbf{Z}$). For $k\in \textbf{Z}$, define
 $$h_k = \textrm{sign}(g-g^{\nu_k}),\qquad a^k= 2^k(h_k-h_k^{\nu_k}).$$
 It is easy to see that each $a^k(k\in \textbf{Z})$ is a w-3-atom.
 Thus, by Theorem 3.3, if $f^N$ is again defined by (5.1), then $\|f^N\|_ {wH_\Phi}\leq C.$
 Therefore
 \be
 C\|l\|&\geq &|l(f^N)|=|\E(f^N g)|=|\sum_{k=-N}^N \E(a^k g)|=|\sum_{k=-N}^N 2^k
 \E((h_k-h_k^{\nu_k})g)|\\
 & = & |\sum_{k=-N}^N 2^k \E(h_k(g-g^{\nu_k}))|= \sum_{k=-N}^N 2^k
 \|g-g^{\nu_k}\|_1\\
 & = &\sum_{k=-N}^N
 \frac{1}{\phi\big(\frac{1}{\Phi(2^k)}\big)}
 \Big(\frac{1}{\Phi(2^k)}\Big)^{-1}
 \|g-g^{\nu_k}\|_1
 \ee
 Taking over all $N\in \textbf{N}$ and the supremum over all of such
 stopping times such that $P(\nu_k <\infty)\leq (\Phi(2^k))^{-1}$, $k\in
 \textbf{Z}$, we obtain
 $$\|g\|_{w\mathcal {L}_{1,\phi}}=\int_{0}^{\infty}\frac{t_{\phi}^{1}(x)}{x}dx \leq C\sum_{k=-\infty}^{\infty}t_{\phi}^{1}\Big(\frac{1}{\Phi(2^k)}\Big)\leq C\|l\|.$$
 The proof of Lemma 5.2 is complete.

 \begin{lem}
 Let $\Phi \in \mathcal {G}_\ell$ with $\ell
 \in (0,1]$, $q_{\Phi^{-1}}<\infty$ and
 $\phi(r)=1/(r\Phi^{-1}(1/r))$. Then for $q\in[1,\infty)$, if $\{\mathcal {F}_n\}_{n\geq 0}$ is
 regular, we have $$(w\mathscr{H}_\Phi)^*= w\mathcal {L}_{q,\phi}.$$
 \end{lem}

 \noindent {\bf Proof.} If $g\in w\mathcal {L}_{q,\phi} $ and $$l_g(f):=\E(fg),\qquad f\in
 L_{q'},$$ where $q'=q/(q-1)$, then by Proposition 2.10 and Theorem
 5.2 we have
 $$|l_g (f)|=|\E(fg)| \leq C \|f\|_{wH_\Phi} \|g\|_{w\mathcal {L}_{1,\phi}} \leq C \|f\|_{wH_\Phi} \|g\|_{w\mathcal {L}_{q,\phi}}.$$

 Conversely, suppose that $l\in(w\mathscr{H}_\Phi)^* $. Since $L_{q'}\subset L_1 \subset L_\Phi \subset
 w\mathscr{L}_\Phi$,  $L_{q'}$ can be embedded continuously in $w\mathscr{H}_\Phi$.
 Thus there exists $g\in L_q$ such that $l$ equals $l_g$ on $L_{q'}$.
 Let $\nu_k$ be the stopping times satisfying $P(\nu_k <\infty)\leq
 2^{-k}$ ($k\in \textbf{Z}$). Define
 $$h_k= \frac{|g-g^{\nu_k}|^{q-1} \textrm{sign} (g-g^{\nu_k}) }  {\|g-g^{\nu_k}\|_q^{q-1}}, \qquad a^k=\Phi^{-1}(2^k) (2^k)^{-1/q'}(h_k - h_k^{\nu_k}).$$
 Then $\|h_k\|_{q'}=1$ and $a^k =0$ on the set $\{\nu_k=\infty\}$. For $\lambda>0$, choose $j\in \textbf{Z}$
 such that $2^{j}\leq \lambda < 2^{j+1}$.
 Define $f^N$, $g^N$ and $h^N$ ($N\in \textbf{N}$) again by (5.1),
 then
 $$\|M(g^N)\|_{q'} \leq \|g^N\|_{q'} \leq \sum_{k=-N}^{j-1}\|a^k\|_{q'} \leq 2\sum_{k=-N}^{j-1}\Phi^{-1}(2^k)(2^k)^{-1/q'}$$
 and
 \be
 P(M(g^N)> \Phi^{-1}(\lambda)) & \leq & \frac
 {1}{(\Phi^{-1}(\lambda))^{q'}} \|M(g^N)\|_{q'}^{q'}\\
 & \leq & \frac{2^{q'}}{(\Phi^{-1}(\lambda))^{q'}}
 \Bigg(\sum_{k=-N}^{j-1}\Phi^{-1}(2^k)(2^k)^{-1/q'}\Bigg)^{q'}\\
 & \leq & 2^{q'} \Bigg(\sum_{k=-N}^{j-1}\frac{(2^{-k})^{1/q'}\Phi^{-1}(2^k)}{\Phi^{-1}(\lambda)}\Bigg)^{q'}\\
 & \leq & C\Bigg(\sum_{k=-N}^{j-1}(2^{-k})^{1/q'}\Big(\frac{2^k}{\lambda}\Big)^p\Bigg)^{q'}\\
 &=& C\cdot\lambda^{-q'p}\Big(\sum_{k=-N}^{j-1}(2^{p-\frac{1}{q'}})^k\Big)^2\leq
 C\lambda^{-1}.
 \ee
 The last inequality holds since $1/q'<1 \leq p$.
 Applying the method used in Theorem 5.1, we conclude that
 $\|f^N\|_{wH_\Phi}\leq C$.
 Consequently,
 \be
 C\|l\|&\geq &|l(f^N)|=|\sum_{k=-N}^N \E(a^k g)|=|\sum_{k=-N}^N
 \Phi^{-1}(2^k)(2^k)^{-1/q'}
 \E((h_k-h_k^{\nu_k})g)|\\
 & = & |\sum_{k=-N}^N \Phi^{-1}(2^k) (2^k)^{-1/q'}\E(h_k(g-g^{\nu_k}))|= \sum_{k=-N}^N \Phi^{-1}(2^k)(2^k)^{-1/q'}
 \|g-g^{\nu_k}\|_q\\
 & = &\sum_{k=-N}^N \frac{1}{\phi(2^{-k})}(2^{-k})^{-1/q}
 \|g-g^{\nu_k}\|_q
 \ee
 Taking over all $N\in \textbf{N}$ and the supremum over all of such
 stopping times such that $P(\nu_k <\infty)\leq 2^{-k}$($k\in
 \textbf{Z}$), we obtain
 $$\|g\|_{w\mathcal {L}_{q,\phi}}=\int_{0}^{\infty}\frac{t_{\phi}^{q}(x)}{x}dx \leq C\sum_{k=-\infty}^{\infty}t_{\phi}^{q}(2^{-k})\leq C\|l\|.$$
 The proof of the theorem is complete.\\

We finally formulate the weak version of the John-Nirenberg theorem,
which directly results from Lemma 5.2 and Lemma 5.3.
 \begin{thm}
 If there exists $\Phi\in \mathcal {G}_\ell$ with $q_{\Phi^{-1}}<\infty$ such that $\phi(r)= \frac{1}{r\Phi^{-1}(1/r)}$ for all $r\in(0,\infty)$, and $\{\mathcal {F}_n\}_{n\geq 0}$ is
 regular.
 Then $w\mathcal {L}_{q,\phi}$ spaces are equivalent for all $1\leq q<\infty$.
 \end{thm}

 \begin{remark} Considering $\Phi(t)\equiv1,$ we obtain the
 John-Nirenberg theorem, Corollary 8 in \cite{weisz1} due to Weisz.
 \end{remark}

\end{document}